\documentclass{amsart}

\usepackage{amsmath,amsfonts,amssymb}

\newcommand \chr{{\mathrm {char}}}

\begin{document}

\title[Orthogonal Lie algebras in dimension 4]{A note on orthogonal Lie algebras in dimension 4\\viewed as current Lie algebras}

\author{Martin Chaktoura}
\address{Department of Mathematics and Statistics, Univeristy of Regina, Canada}
\email{martin\_chaktoura@yahoo.com.ar}

\author{Fernando Szechtman}
\address{Department of Mathematics and Statistics, Univeristy of Regina, Canada}
\email{fernando.szechtman@gmail.com}
\thanks{The second author was supported in part by an NSERC discovery grant}




%
%



\begin{abstract}
Orthogonal Lie algebras in dimension 4 are identified as current
Lie algebras, thus producing a natural decomposition for them over
any field.
\end{abstract}

\maketitle

Let $f:V\times V\to F$ be a non-degenerate symmetric bilinear form
defined on a vector space $V$ of dimension 4 over a field $F$. If
$\chr(F)=2$ we further assume that $f$ is not alternating. Let $D$
be the discriminant of $f$ relative to a basis of~$V$. Let $L(f)$
be the subalgebra of $\mathfrak{gl}(V)$ associated to $f$, formed
by all $x\in\mathfrak{gl}(V)$ that are skew-adjoint relative
to~$f$. Let $M=[L(f),L(f)]$, which coincides with $L(f)$ if and
only if $\chr(F)\neq 2$. In any case, $M$ is a 6-dimensional
orthogonal Lie algebra.

According to \cite{B}, Chapter 1, \S 6, Exercise 26, if
$\chr(F)\neq 2$ then $M$ is either the direct sum of two
3-dimensional simple ideals or simple, depending on whether $D$ is
a square in $F$ or not. On the other hand, if $\chr(F)=2$ then
either $M=N\ltimes R$ or $M$ is simple, depending, again, on
whether $D$ is a square in $F$ or not; here $N$ is a simple
3-dimensional simple subalgebra and $R$ is the solvable radical of
$M$.

We may view these orthogonal Lie algebras as current Lie algebras,
and in this way explain all cases described above in a uniform
manner.

\medskip

\noindent{\bf Theorem. } We have $M\cong [L(f|_W),L(f|_W)]\otimes
F[X]/(X^2-D)$, where $W$ is an arbitrary 3-dimensional subspace of
$V$ such that $f|_W$ is non-degenerate.

\begin{proof}
It follows from \cite{K}, Theorems 4 and
20, that $V$ admits an orthogonal basis $B=\{v_1,v_2,v_3,v_4\}$.
Thus, the Gram matrix of $f$ relative to $B$ is diagonal with
non-zero entries $a,b,c,d$ and $D=abcd$. Let
$$
f_1=be_{12}-ae_{21},\; f_2=ce_{23}-be_{32},\; f_3=ce_{13}-ae_{31},
$$
$$
h_1=ab(de_{34}-ce_{43}),\; h_2=bc(de_{14}-ae_{41}),\;
h_3=ac(be_{42}-de_{24}).
$$
Then $f_1,f_2,f_3,h_1,h_2,h_3$ is a basis of $M$, with
multiplication table
$$
[f_1,f_2]=bf_3,\; [f_2,f_3]=cf_1,\; [f_3,f_1]=af_2,
$$
$$
[f_1,h_2]=bh_3,\; [f_2,h_3]=ch_1,\; [f_3,h_1]=ah_2,
$$
$$
[f_2,h_1]=-bh_3,\; [f_3,h_2]=-ch_1,\; [f_1,h_3]=-ah_2,
$$
and
$$
[h_1,h_2]=Dbf_3,\; [h_2,h_3]=Dcf_1,\; [h_3,h_1]=Daf_2.
$$
Thus $M$ has the same multiplication table as
$[L(f|_W),L(f|_W)]\otimes F[X]/(X^2-D)$, where $W$ is the span by
$v_1,v_2,v_3$. By \cite{K}, Theorems 4 and 20, $W$ can be replaced
by any 3-dimensional subspace of $V$ where the restriction of $f$
is non-degenerate.
\end{proof}

\medskip

It is obvious from the Theorem that $M$ decomposes exactly as
prescribed above if~$D$ is a square. Suppose $D$ is not a square.
Then $K=F[X]/(X^2-D)$ is a quadratic field extension of $F$. Since
$[L(f|_W),L(f|_W)]$ is a perfect 3-dimensional Lie algebra over
$F$, it follows that $[L(f|_W),L(f|_W)]\otimes K$ is a simple Lie
algebra over $K$, and hence over $F$, as seen below. Thus, by the
Theorem, $M$ is a simple Lie algebra over $F$.

\medskip

\noindent{\bf Lemma. }(cf. \cite{LP}, Lemma 2.7) Let $L$ be a
simple Lie algebra over a field $K$. Suppose that $F$ is a
subfield of $K$. Then $L$ is simple over $F$.

\begin{proof}
Let $I$ be a non-zero ideal of $L$ over
$F$. The $K$-span of $I$ is a non-zero ideal $J$ of $L$ over $K$,
so $J=L=[L,L]=[L,J]=[L,I]\subseteq I$.
\end{proof}

\medskip

A simple Lie algebra need not remain simple, or even semisimple,
upon field extension. For instance, let $L$ be an absolutely
simple Lie algebra over an imperfect field $F$ of characteristic
$p$, and let $s\in F\setminus F^p$ and $K=F[X]/(X^p-s)$, which has
an element $t$ satisfying $t^p=s$. Then $P=L\otimes K$ is simple
over $K$ and hence over~$F$, but $P\otimes K\cong L\otimes
K[X]/(X^p-s)\cong L\otimes K[X]/(X-t)^p$ is not semisimple over
$K$. 

\medskip

\noindent{\bf Acknowledgment.} We are grateful to A. Pianzola for
pointing out reference \cite{LP}, which replaced a more convoluted
argument.



\begin{thebibliography}{99}

\bibitem[B]{B}
N. Bourbaki,
\newblock{\em Lie groups and Lie algebras, Chapters 1-3},
\newblock Springer-Verlag, Berlin, 1989.

\bibitem[K]{K}
I. Kaplansky,
\newblock{\em Linear algebra and geometry.
A second course.},
\newblock Allyn and Bacon, Boston, 1969.



\bibitem[LP]{LP}
M. Lau and A. Pianzola,
\newblock{\em Maximal Ideals and Representations of Twisted Forms
of Algebras},
to appear in Algebra and Number Theory.

\end{thebibliography}
\end{document}